\documentclass{rmmcart}

\usepackage{lineno,hyperref}
\usepackage{amsmath, amssymb, amscd, amsthm}
\usepackage{accents}
\usepackage{setspace}
\usepackage{graphicx}
\usepackage{subfigure}
\usepackage{psfrag}
\usepackage{epstopdf}
\usepackage{textcomp}
\usepackage{xspace}
\usepackage{booktabs}
\usepackage{multirow}
\usepackage{tabu}
\usepackage{enumerate}
\usepackage{tikz}
\usepackage{tikz-cd}
\usepackage{xcolor}
\usepackage{color}

\newtheorem{thm}{Theorem}[section]

\newtheorem{cor}[thm]{Corollary}
\newtheorem{lem}[thm]{Lemma}
\newtheorem{prop}[thm]{Proposition}

\newtheorem{defn}[thm]{Definition}

\newtheorem{exam}[thm]{Example}

\newtheorem*{prop*}{Proposition}
\newtheorem*{cor*}{Corollary}
\newtheorem*{thm*}{Theorem}
\newtheorem*{defn*}{Definition}

\def\pf{\bigskip\noindent {\bf {Proof}  }~}

\def\lemma{\noindent {\bf{Lemma} }~}

\def\qed{ \begin{flushright} $\square$ \end{flushright} }
\def\bf{\textbf}

\newcommand{\Hom}{\text{Hom}}
\newcommand{\coker}{\text{coker }}

\newcommand{\im}{\text{im }}

\newcommand{\Tor}{\text{Tor}}
\newcommand{\Ext}{\text{Ext}}
\newcommand{\depth}{\text{depth}}

\newcommand{\ann}{\text{Ann}}
\newcommand{\g}{\text{grade}}

\newcommand{\Epi}{\text{Epi}}
\newcommand{\Supp}{\text{Supp}}
\newcommand{\G}{\text{G}}

\title{Linkage and Intermediate \(C\)-Gorenstein Dimensions}

\author{Joseph P.\ Brennan}
\email{Joseph.Brennan@ucf.edu}
\address{Department of Mathematics\\ University of Central Florida\\  Orlando, FL 32816\\ USA}

\author{Alexander York}  
\email{Alexander.York@ucf.edu}
\address{Department of Mathematics\\ University of Central Florida\\  Orlando, FL 32816\\ USA}
\thanks{York Current Email: A.York@ufl.edu}


\begin{document}
\begin{abstract}
This paper brings together two theories in algebra that have had been extensively developed in recent years.  First is the study of various homological dimensions and what information such invariants can give about a ring and its modules.  A collection of intermediate \(C\)-Gorenstein dimensions is defined and this allows generalizations of results concerning \(C\)-Gorenstein dimension and certain Serre-like conditions.  Second is the theory of linkage first introduced by Peskine and Szpiro and generalized to modules by Martinskovsky and Strooker.  Using the further generalization of module linkage of Nagel, results are proven connecting linkage with these homological dimensions and Serre-like conditions. 
\end{abstract}
\maketitle

\setcounter{section}{0}


\section{Introduction}\label{Introduction}

The theory of linkage of ideals was first studied by Peskine and Szpiro  \cite{LiaisonPeskineSzpiro}. There is a voluminous literature concerning linkage of ideals 
\cite{LiaisonPeskineSzpiro,LinkageMartin,NotesLiaisonSchenzel,StructureLinkageHuneke,AlgebraicLinkageHuneke,PureLinkageKustinMiller,LiaisonHartshorne}.  
Complete intersection and Gorenstein ideal linkage was generalized by Martin \cite{LinkageMartin} to module linkage and extended by Martinskovsky and Strooker \cite{LinkageMartsinkovskyStrooker} using projective presentations.  
Martinskovsky and Strooker take a projective presentation of an \(R\)-module \(M\)
\[
P_1 \longrightarrow P_0 \longrightarrow M \longrightarrow 0
\]
and consider the module \(\lambda M=\coker\!(\Hom_R(M,R)\to\Hom_R(P_0,R))\).  Then two modules \(M\) and \(N\) are directly linked if \(M\cong \lambda N\) and \(N\cong \lambda M.\)  This generalizes complete intersection linkage for ideals \cite[Proposition 1]{LinkageMartsinkovskyStrooker}.\\

The linkage introduced by Nagel \cite{NagelLiaison} generalized Gorenstein linkage for ideals even further.  This is accomplished by using modules that play a role that is analogous to that played by Gorenstein ideals. Such modules are called quasi-Gorenstein modules.  An \(R\)-module \(Q\) of grade \(q\) is quasi-Gorenstein if the following hold:
\begin{itemize}
\item[\((i)\)] \(\Ext_R^{q+i}(Q,R)=0=\Ext_R^{q+i}(\Ext_R^{q}(Q,R),R)\) for \(i>0\)
\item[\((ii)\)] The natural map \(Q\to\Ext_R^q(\Ext_R^q(Q,R),R)\) is an isomorphism
\item[$(iii)$] There is some isomorphism \(\alpha:Q\to\Ext_R^q(Q,R).\)
\end{itemize}
One can then take a short exact sequence
\[
\begin{tikzcd}
0 \arrow[r] & K \arrow[r, "\psi"] & Q \arrow[r, "\varphi"] & M \arrow[r] & 0
\end{tikzcd}
\]
of $R$-modules where $Q$ is quasi-Gorenstein, \(q=\g_R(Q)=\g_R(M).\)  Then, using the isomorphism \(\alpha\colon Q\to\Ext_R^q(Q,R),\) we obtain the following short exact sequence after dualizing
\[
\begin{tikzcd}
0 \arrow[r] & \Ext_R^q(M,R) \arrow[r] & Q \arrow[r] & \mathcal{L}_Q(M) \arrow[r] & 0
\end{tikzcd}
\]
where we write \(\mathcal{L}_Q(M)=\im\!((\psi)^*\circ\alpha)\) which is a submodule of \(\Ext_R^q(K,R).\)  Then one says that two modules \(M\) and \(N\) are directly linked by \(Q\) if \(\mathcal{L}_Q(M)\cong N\) and \(\mathcal{L}_Q(N)\cong M.\) Nagel showed that many of the results concerning Gorenstein ideal linkage generalize to this setting.  Our goal in this paper is to show that results following from other types of module linkage also generalize to this setting, allowing for a more unified theory of module linkage.  Further, we will move to a linkage involving a more general collection of modules than that of quasi-Gorenstein modules, which will be called \(C\)-quasi-Gorenstein modules where \(C\) is a semidualizing \(R\)-module (Definition \ref{semidualizingdefn}).\\

One of the themes that is present in this construction is that we choose to ``represent" our module by a \(C\)-quasi-Gorenstein module of the same grade.  As we shall see, this illuminates the intuition and motivation concerning previous results and helps us wade through the difficulties that can be incident to comparing modules and ideals in different codimensions.  Before addressing the issue of linkage, we first construct the tools that are necessary to state our results.  This is done through the study of homological dimensions and their use in considering Serre-like conditions.\newline

The study of homological dimensions can be traced back to the earliest days of homological algebra but the first treatise that impacts this paper is the celebrated ``Stable Module Theory" of Auslander and Bridger \cite{AuslanderBridgerSMT}.  They build upon the ideas of projective, flat, and injective dimension by defining \(G\)-dimension to characterize Gorenstein rings.  Recall that an \(R\)-module \(M\) is \textit{totally reflexive} or has \textit{Gorenstein dimension zero} if the following hold:
\begin{itemize}
\item[\((a)\)] \(M\) is finitely generated
\item[\((b)\)] The natural biduality map \(M\to\Hom_R(\Hom_R(M,R),R)\) is an isomorphism
\item[\((c)\)] \(\Ext_R^{i}(M,R)=0=\Ext_R^{i}(\Hom_R(M,R),R)\) for \(i>0\)
\end{itemize} 
One can then define \(G\)-dimension of any \(R\)-module by resolving the module using totally reflexive modules.  There has been considerable focus on generalizing such a construction to characterize other types of rings \cite{CMDIMChristensenFoxbyFrankid, CMDIMHolmJorgensen, CIDIM, GIPFDIMEnochsJenda, GDimIyengar}.  However, in this paper we follow a different train of thought.\newline

Over a Noetherian ring \(R,\) an \(R\)-module \(C\) is said to be \textit{semidualizing} if the natural map \(R\to\Hom_R(C,C)\) is an isomorphism and \(\Ext_R^{i}(C,C)=0\) for \(i\geq 1.\)  One can then ``replace" \(R\) with \(C\) (in the appropriate places) in the definition of Gorenstein dimension zero to define \(C\)-Gorenstein dimension zero.  This was first explored by Foxby \cite{QuasiPerfectModFoxby} and Golod \cite{GDimandGenPerIGolod}.  Many results generalize nicely to this new dimension, but as briefly noted above, every \(C\)-Gorenstein dimension zero module has grade zero and it would help if we could choose what grade they have.  Doing this will allow us to choose the grade in which our \(C\)-Gorenstein resolutions will sit.  We will give a short but illuminating example (Example \ref{GCiExample}) in the next section motivating our construction.  We will define an \(R\)-module, \(M,\) to have \(G_C^j\)-dimension zero as follows:
\begin{defn*}[\textbf{\ref{GCidimzero}}]
An \(R\)-module \(M\) of grade \(j\) is said to have \(G_C^j\)-dimension zero if the following hold:
\begin{itemize}
\item[\((i)\)] \(M\) is fintely generated
\item[\((ii)\)] The natural map \(M\to\Ext_R^j(\Ext_R^j(M,C),C)\) is an isomorphism
\item[\((iii)\)] \(\Ext_R^i(M,C)=0\) for \(i> j\)
\item[\((iv)\)] \(\Ext_R^i(\Ext_R^j(M,C),C)=0\) for \(i\neq j.\)
\end{itemize}
\end{defn*}

\noindent Using this definition, one can define new homological dimensions just as was done for \(C\)-Gorenstein dimension.  Many results about \(C\)-Gorenstein dimension generalize in exactly the expected manner.  We then use these \(G_C^j\)-dimensions to prove results analogous to those in \cite{LinkageModulesSerreSadDib} with regards to generalized Serre conditions and linkage. \newline

We will now discuss the organization of the paper.  In Section 2, we define the intermediate $C$-Gorenstein dimensions  by resolving modules using $G_C^j$-dimension zero modules.  Then we prove results generalizing $C$-Gorenstein dimension.  Specifically, we prove the following
\begin{thm*}[\textbf{\ref{GCjdimthm}}]
If an \(R\)-module \(M\) has finite \(G_C^j\)-dimension then \(G_C^j-\hbox{\rm dim}_R(M)=\alpha_R(M)-j,\) where \(\alpha_R(M)=\sup\{i~|~\Ext_R^i(M,R)\neq 0\}.\)
\end{thm*}
From this one obtains a chain of inequalities
\[
G_C^{\g_C(M)}\text{-dim}_R(M)\leq G_C^{\g_C(M)-1}\text{-dim}_R(M)\leq \cdots \leq G_C^{1}\text{-dim}_R(M) \leq G_C^{0}\text{-dim}_R(M).
\]
It is natural to ask if these dimensions are simultaneously finite, and we provide that result in Theorem \ref{GCdimallalwaysfinite}.  Finally, we state an Auslander-Bridger type formula for these dimensions as follows
\begin{cor*}[\textbf{\ref{GCjdimdepthequation}}]
If \((R,\mathfrak{m},k)\) is a local ring and \(M\) has finite \(G_C^j\)-dimension, then
\[
G_C^j\text{-dim}_R(M)=\depth(R)-\depth_R(M)+j.
\]
\end{cor*}

In Section 3, we follow the ideas of Sadeghi and Dibaei \cite{LinkageModulesSerreSadDib} and focus on the connection between certain Serre-like conditions and these homological dimensions.  We define a modified Serre-type condition \(\tilde{S}_n^g\) as follows.  An \(R\)-module \(M\) satisfies \(\tilde{S}_n^g\) if \(\depth_{R_p}(M_p)+g\geq \min\{n,\depth(R_p)\}\) for all \(p\in\text{Spec}(R).\)  If \(g=0\) this is the condition \(\tilde{S}_n.\)  We then prove results analogous to those in \cite{LinkageModulesSerreSadDib}.  The main theorem of this section is the following
\begin{thm*}[\textbf{\ref{GCSerreconditionthm}}]
Let \(C\) be a semidualizing \(R\)-module, \(M\) an \(R\)-module with \(\g_R(M)=g,\) \(n\geq g,\) and \(M\) having locally finite \(G^g_C\)-dimension. Then the following are equivalent:
\begin{itemize}
\item[\((i)\)] \(\Ext_R^{g+i}(D_C^gM,C)=0\) for \(1\leq i\leq n-g\)
\item[\((ii)\)] \(M\) is an \((n-g)^{\text{th}}\) \(\mathcal{G}_C^g\)-syzygy
\item[\((iii)\)] \(M\) satisfies \(\tilde{S}_n^g\)
\item[\((iv)\)] \(\g_{C_p}(\Ext_{R_p}^{\g_{C_p}(M_p)+i}(M_p,C_p))\geq i+n\) for \(i\geq 1\) and \(p\in\text{Spec}(R)\) where \(\depth(R_p)\leq i+n-1.\)
\end{itemize}
\end{thm*}
\noindent where \(D_C^gM\) is an Auslander-type dual to \(M\) in grade \(g.\)  This result corresponds to \cite[Proposition 2.4]{LinkageModulesSerreSadDib}. We end this section by proving that for any \(R\)-module, \(M,\) in the Auslander class of \(C\) with locally finite \(G_C^g\)-dimension
\[
\Ext_R^{g+i}(D_R^gM,R)=0\text{ for }1\leq i\leq n-g\quad\Leftrightarrow\quad \Ext_R^{g+i}(D_C^gM,C)=0\text{ for }1\leq i\leq n-g
\]
generalizing \cite[Theorem 2.12]{LinkageModulesSerreSadDib}.\newline

In Section 4, we meld together the techniques and ideas from \cite{LinkageModulesSerreSadDib}, \cite{LinkageMartsinkovskyStrooker}, and \cite{NagelLiaison} to connect together the theory of module linkage with these homological dimensions and Serre-like conditions.  As was stated above, Nagel generalized module linkage by allowing non-projective modules to present \(M.\)  In particular, he defines module linkage using quasi-Gorenstein modules, from which we use \(C\)-quasi-Gorenstein modules.  Then given a surjective map \(Q\to M\to 0\) where \(Q\) is \(C\)-quasi-Gorenstein with \(\g_R(Q)=\g_R(M)\) one can define \(\mathcal{L}_Q(M)\) as we stated earlier.  With this we prove a corresponding result to \cite[Theorem 2]{LinkageMartsinkovskyStrooker}
\begin{thm*}[\textbf{\ref{horizontallylinkedtheorem}}]
An \(R\)-module \(M\) of grade \(g\) is horizontally linked by a \(C\)-quasi-Gorenstein \(R\)-module \(Q\) if and only if \(M\) is \(\mathcal{G}_C^g\)-stable and \(\Ext_R^{g+1}(D_C^gM,C)=0.\)
\end{thm*} 
We then prove Proposition \ref{LinkageandSerreLocalDuality1} and Corollary \ref{LinkageandSerreLocalDuality2} which are generalizations of \cite[Proposition 2.6]{LinkageModulesSerreSadDib} and \cite[Corollary 2.8]{LinkageModulesSerreSadDib}, which are themselves generalizations of a result of Schenzel \cite[Theorem 4.1]{NotesLiaisonSchenzel}. The paper culminates in Theorem \ref{GCDimandLinkageThm}, a generalization of \cite[Corollary 2.14]{LinkageModulesSerreSadDib}, which can be stated as follows:
\begin{thm*}[\textbf{\ref{GCDimandLinkageThm}}]
Supppose that \(R\) is a Cohen-Macaulay ring \(R\) with semidualizing module \(C\) and  \(M\) is a \(\mathcal{G}_C^g\)-stable \(R\)-module of grade \(g.\)  Let \(M\) be in the Auslander class of $C$ and have locally finite \(G_C^g\)-dimension then the following are equivalent:
\begin{itemize}
\item[\((i)\)] \(M\) satsifies \(\tilde{S}_n^g\) 
\item[\((ii)\)] \(M\) is horizontally linked by some \(C\)-quasi-Gorenstein \(R\)-module \(Q\) and \(\Ext_R^{g+i}(\mathcal{L}_Q(M),C)=0\) for \(0<i<n-g.\)
\end{itemize}
\end{thm*}


\section{Intermediate \(C\)-Gorenstein Dimensions}\label{ICGD}

Throughout this section \(R\) will be a Noetherian ring and all \(R\)-modules are assumed to be finitely generated.  Definitions and results stated without reference or proof can be found in \cite{EisenbudCA,GDIMChristensen,BasicAlgebra1,BasicAlgebra2,QuasiPerfectModFoxby}.  The goal of this section is to build up a theory of dimensions for modules of different grades which will allow us to generalize results concerning generalized Serre conditions and module linkage.  Since we feel that the results concerning these ideas have their own merit we have presented them in this section.\newline

First, we recall the definition of a semidualizing \(R\)-module, which were first studied in \cite{QuasiPerfectModFoxby} and \cite{GDimandGenPerIGolod}
\begin{defn}\label{semidualizingdefn}
An \(R\)-module \(C\) is called \textbf{semidualizing} if the homothety map \(R \to \Hom_R(C,C)\) is an isomorphism and \(\Ext_R^i(C,C)=0\) for \(i>0.\)
\end{defn}

Obvious examples of semidualizing modules are the ring \(R\) and the canonical modules \(\omega_R\) of a Cohen-Macaulay ring \(R.\)  There has been a large body of research into semidualizing modules over rings and their connections with homolgical dimensions, \cite{RelativeSemidualizingSather-Wagstaff, SemidualizingSather-Wagstaff, QuasiPerfectModFoxby, GDimandGenPerIGolod}.\newline

For a semidualizing \(R\)-module \(C,\) we will call \[\depth_R(\ann_R(M),C)=\inf\{i~|~\Ext_R^i(M,C)\neq 0\}\] the \(C\)-grade of \(M,\) denoted by \(\g_C(M).\)  Note that \(\g_R(M)=\g_C(M)\) for any semidualizing 
\(R\)-module \(C,\) and so we will stick with \(\g_R(M).\)  Also we will use \[\alpha_C(M)=\sup\{i~|~\Ext_R^i(M,C)\neq 0\}.\]  Once again, \(\alpha_C(M)=\alpha_R(M)\) for any semidualizing \(R\)-module \(C\) and so we will use \(\alpha_R(M).\)  Given an \(R\)-module \(M\) with \(\g_R(M)\geq g,\) it follows from results in \cite{AuslanderBridgerSMT} or \cite{FossumDGR} that there is a natural map
\[
\begin{tikzcd}
\delta_C^g(M): M \arrow[r] & \Ext_R^g(\Ext_R^g(M,C),C)
\end{tikzcd}
\]
we denote by \(\delta_C^g(M).\)\newline

The Gorenstein dimension of a module has been extended to the \(C\)-Gorenstein dimension of a module \cite{GorensteinModandRelatedFoxby}.  One of the disadvantages of \(C\)-Gorenstein dimension is that the modules that comprise the resolutions sit in grade 0 whereas the module being resolved may not.  It is useful to be able to use modules in the same grade as it helps sift out unwanted information and makes proof techniques simpler.  We can illustrate this with an example.  
\begin{exam}\label{GCiExample}
Consider the ring \(R=k[x,y]\) where \(k\) is an infinite field and let \(C=R\) be the semidualizing module and \(k\) be the \(R\)-module under consideration.  We can take a \(G_C\)-resolution of \(k\) as
\[
\begin{tikzcd}[ampersand replacement = \&]
0 \arrow[r] \& R \arrow{r}{\begin{pmatrix} x \\ y \end{pmatrix}} \& R^2 \arrow{r}{\begin{pmatrix} -y & x\end{pmatrix}} \& R \arrow[r] \& k \arrow[r] \& 0
\end{tikzcd}
\]
which is the Koszul complex of \(k.\)  If we also consider the \(R\)-module \(R/(x),\) we see that these modules sit in different grades, \(k\) in grade 2 and \(R/(x)\) in grade 1.  If we wanted to consider the \(R\)-module \(\Hom_R(k, R/(x)),\) as one often does, we can apply \(\Hom_R(-,R/(x))\) to the Koszul complex  of \(k.\)  However, it may be simpler instead to consider the following way of representing \(k:\)
\[
\begin{tikzcd}
0 \arrow[r] & R/(x) \arrow[r, "\cdot y"] & R/(x) \arrow[r] & k \arrow[r] & 0.
\end{tikzcd}
\]
Notice that since \(R/(x)\) is of grade 1, both modules presenting \(k\) are of grade 1.  Then it is much simpler to understand \(\Hom_R(k,R/(x))\) by using this short exact sequence.  In fact, we have
\[
\begin{tikzcd}
0 \arrow[r] & \Hom_R(k,R/(x)) \arrow[r] & \Hom_R(R/(x),R/(x)) \arrow[r, "\cdot y"] & \Hom_R(R/(x),R/(x))
\end{tikzcd}
\]
which shows that \(\Hom_R(k,R/(x))\) is the kernel of the right map.  This shows that there is an advantage to using representations of modules in higher grades. 

Notice that \(R/(x)\) as an \(R\)-module satisfies
\begin{itemize}
\item[\((i)\)] \(\g_R(R/(x))=1\)
\item[\((ii)\)] \(\Hom_R(R/(x),R)=0,\) \(\Hom_R(\Ext_R^1(R/(x),R),R)=0,\) and \(\Ext_R^2(R/(x),R)=0\)
\item[\((iii)\)] \(R/(x)\cong \Ext_R^1(\Ext_R^1(R/(x),R),R).\)
\end{itemize}   
In other words, these conditions are similar to those of \(G_C\)-dimension zero but for a higher index in Ext.
\end{exam}  

This motivates the following definition.

\begin{defn}\label{GCidimzero}
Suppose that \(C\) is a semidualizing \(R\)-module.  We say that an \(R\)-module \(M\) with \(\g_R(M)=g\) has \(\mathbf{\G_C^g}\)\textbf{-dimension zero}, or \(\G_C^g\)-dim\(_R(M)=0,\) if the following conditions hold:
\begin{itemize}
\item[\((i)\)] \(\Ext_R^{g+i}(M,C)=0\) for \(i>0\)
\item[\((ii)\)] \(\Ext_R^{g+i}(\Ext_R^g(M,C),C)=0\) for \(i> 0\)
\item[\((iii)\)] The map \(\delta_C^g(M):M\to \Ext_R^g(\Ext_R^g(M,C),C)\) is an isomorphism
\end{itemize}
\end{defn}

We will call the class of all \(\G_C^g\)-dimension zero modules \(\mathcal{G}_C^g.\)  Notice that if \(M\in\mathcal{G}_C^j\) then so is \(\Ext_R^{j}(M,C).\)  Then just as with \(C\)-Gorenstein dimension we can resolve any module with grade greater than or equal to \(g\) by modules in \(\mathcal{G}_C^g\) and define

\begin{defn}\label{GCidim}
Suppose that \(C\) is a semidualizing \(R\)-module and \(\g_R(M)\geq j\) for an \(R\)-module \(M.\)  Given an exact sequence
\[
\begin{tikzcd}
\mathcal{C}: 0 \arrow[r] & M_n \arrow[r] & M_{n-1} \arrow[r] & \cdots \arrow[r] & M_1 \arrow[r] & M_0 \arrow[r] & M \arrow[r] & 0
\end{tikzcd}
\]
where \(M_i\in \mathcal{G}_C^j\) for \(i=0,\ldots, n,\) we will say that \(\mathcal{C}\) is a \(\G_C^j\)-resolution of \(M.\)
\end{defn}

If \(M\) has a \(\G_C^j\)-resolution of length \(n+1,\) such as the one above, then we say that the \textbf{jth intermediate} \(\mathbf{C-}\)\textbf{Gorenstein dimension} of \(M\) is less than or equal to \(n,\) 
or \(\G_C^j\)-dim\(_R(M)\leq n.\)  Note that this agrees exactly with \(\G_C\)-dim\(_R(M)\) when \(j=0.\)  A question one may have after this definition is if every module has such a resolution.  
The assumption that \(\g_R(M)\geq j\) was made to provide an affirmative answer to this inquiry.  If \(\g_R(M)\geq j,\) then there exists a \(C\)-regular sequence \(\bar{x}=(x_1,\ldots, x_j)\) in \(\ann_R(M).\)  
Then there is a natural surjective map from \[(R/(\bar{x})\otimes_R C)^\alpha\cong (C/(\bar{x})C)^\alpha \to M\to 0\] for some \(\alpha\geq 0.\)  Clearly \((C/(\bar{x})C)^\alpha\in\mathcal{G}_C^j.\) 
Since \(\g_R(\ker(C/(\bar{x})C\to M)\geq \g_R(M)\geq j\) we continue the process with \(\ker(C/(\bar{x})C\to M).\)  
Continuing in this fashion we can iteratively construct a \(G_C^j\)-resolution of any \(R\)-module \(M\) with \(\g_R(M)\geq j.\)\newline

There are analogous results for \(G_C^j\)-dimension to those on \(G\)-dimension.  We list a few which are necessary for proofs in this paper.  The next result is analogous to \cite[Lemma 1.1.10]{GDIMChristensen}.

\begin{prop}\label{GCShortexactsequence}
Suppose that \(0 \to M' \to M \to M'' \to 0\) is a short exact sequence with \(\g_R(M')=\g_R(M)=\g_R(M'')=g.\)  Then
\begin{itemize}
\item[\((a)\)] If \(M''\in\mathcal{G}_C^g,\) then \(M\in\mathcal{G}_C^g\Leftrightarrow M'\in\mathcal{G}_C^g\)
\item[\((b)\)] If \(M\in\mathcal{G}_C^g,\) then \(\Ext_R^{g+i}(M',C)\cong\Ext_R^{g+i+1}(M'',C)\) for \(i>0\)
\item[\((c)\)] If the sequence splits, then \(M\in\mathcal{G}_C^g\Leftrightarrow M',M''\in\mathcal{G}_C^g.\)
\end{itemize}
\end{prop}  

\pf: Both \((a)\) and \((b)\) are clear by applying \(\Hom_R(-,C)\) and considering the long exact sequence, and \((c)\) follows by the naturality of \(\delta_C^g(-)\) and the commutativity of Ext and direct sum.\qed

The next two lemmas will make the proof of the following results shorter.

\begin{lem}\label{GCExtShiftLem}
Suppose that 
\[
\begin{tikzcd}
0 \arrow[r] & M_n \arrow[r] & M_{n-1} \arrow[r] & \cdots \arrow[r] & M_1 \arrow[r] & M_0 \arrow[r] & M \arrow[r] & 0
\end{tikzcd}
\]
is an exact complex of $R$-module with $M_i\in\mathcal{G}_C^g$ for $i=0,\ldots, n-1$.  Let $K_i=\ker(M_{i-1}\to M_{i-2})$ for $i=2,\ldots, n$ with $K_{n}=M_n$ and $K_{0}=M_0$.  Then
\[
\Ext_R^{\g_R(M_{i-1})+j}(K_i,C)\cong \Ext_R^{\g_R(M_{i-1})+i+j}(M,R)\text{ for }j>0\text{ and }i=0,\ldots, n.
\]
\end{lem}

\pf: Clear by breaking the complex into short exact sequences and applying \(\Hom_R(-,C).\)\qed

\begin{lem}\label{GCAllPerfectFlip}
Suppose that
\[
\begin{tikzcd}
0 \arrow[r] & M_n \arrow[r] & M_{n-1} \arrow[r] & \cdots \arrow[r] & M_1 \arrow[r] & M_0 \arrow[r] & 0
\end{tikzcd}
\]
is exact with \(\g_R(M_i)=\alpha_R(M_i)=j\) for all \(i\geq 0.\)  Then there is an exact sequence
\[
\begin{tikzcd}
0 \arrow[r] & \Ext_R^j(M_0,C) \arrow[r] & \Ext_R^j(M_1,C) \arrow[r] & \cdots \arrow[r] & \Ext_R^j(M_{n-1},C) \arrow[r] & \Ext_R^j(M_n,C) \arrow[r] & 0
\end{tikzcd}
\]
\end{lem}

\pf: Clear by breaking the complex into short exact sequences.\qed

The next result is analogous to \cite[Lemma 1.2.6]{GDIMChristensen}.
\begin{prop}\label{GCFinitedimPerfect0}
Suppose that \(\G_C^j\)-dim\(_R(M)<\infty\) and \(\g_R(M)=j=\alpha_R(M).\)  Then \(M\in\mathcal{G}_C^j.\)
\end{prop}

\pf:  Suppose \(G_C^j\)-dim\(_R(M)=n.\)  We prove by induction on \(n.\)  We are done if \(n=0,\) so suppose \(n=1\) and consider a shortest \(G_C^j\)-resolution of \(M\)
\[
\begin{tikzcd}
0 \arrow[r] & M_1 \arrow[r] & M_0 \arrow[r] & M \arrow[r] & 0.
\end{tikzcd}
\]
By Lemma \ref{GCAllPerfectFlip} we get the following short exact sequence
\[
\begin{tikzcd}
0 \arrow[r] & \Ext_R^j(M,C) \arrow[r] & \Ext_R^j(M_0,C) \arrow[r] & \Ext_R^j(M_1,C) \arrow[r] & 0.
\end{tikzcd}
\]
By Proposition \ref{GCShortexactsequence} \((a)\) we see that \(\Ext_R^j(M,C)\in\mathcal{G}_C^j.\)  Therefore \(\Ext_R^j(\Ext_R^j(M,C),C)\in\mathcal{G}_C^j.\)  
If we apply \(\Hom_R(-,C)\) to the above short exact sequence we get
\[
\begin{tikzcd}
0 \arrow[r] & \Ext_R^j(\Ext_R^j(M_1,C),C) \arrow[r] & \Ext_R^j(\Ext_R^j(M_0,C),C) \arrow[r] & \Ext_R^j(\Ext_R^j(M,C),C) \arrow[r] & 0
\end{tikzcd}
\]
and it is then clear using \(delta_C^j(-)\) that \(M\cong\Ext_R^j(\Ext_R^j(M,C),C)\) and thus \(M\in\mathcal{G}_C^j.\)\newline

Now suppose that \(G_C^j\)-dim\(_R(M)=n\) with \(n>1.\)  Take a shortest \(G_C^j\)-resolution of \(M\)
\[
\begin{tikzcd}
0 \arrow[r] & M_n \arrow[r] & \cdots \arrow[r] & M_0 \arrow[r] & M \arrow[r] & 0
\end{tikzcd}
\]
and let \(K=\ker(M_0\to M).\)  Then \(G_C^j\)-dim\(_R(K)\leq n-1<\infty\) and by applying \(\Hom_R(-,C)\) to the short exact sequence
\[
\begin{tikzcd}
0 \arrow[r] & K \arrow[r] & M_0 \arrow[r] & M \arrow[r] & 0
\end{tikzcd}
\]
we see that \(\Ext_R^{j+i}(K,C)=0\) for \(i>0.\)  Therefore \(\alpha_R(K)=j=\g_R(K)\) and so by the induction hypothesis \(K\in\mathcal{G}_C^j.\)  
Then \(G_C^j\)-dim\(_R(M)\leq 1\) and again by the induction hypothesis \(M\in\mathcal{G}_C^j.\)\qed

Using these results we get the following theorem about \(G_C^g\)-dimension.

\begin{thm}\label{GCjdimthm}
Let \(R\) be a Cohen-Macaulay ring and \(C\) and semidualizing \(R\)-module.  For a finitely generated \(R\)-module \(M\) with \(\g_R(M)\geq j,\) the following are equivalent:
\begin{itemize}
\item[\((i)\)] \(\G_C^j\)-dim\(_R(M)\leq n\)
\item[\((ii)\)] \(\G_C^j\)-dim\(_R(M)<\infty\) and \(n\geq \alpha_R(M)-j\)
\item[\((iii)\)] In any \(\G_C^j\)-resolution,
\[
\begin{tikzcd}
\cdots \arrow[r] & M_i \arrow[r] & M_{i-1} \arrow[r] & \cdots \arrow[r] & M_0 \arrow[r] & M \arrow[r] & 0
\end{tikzcd}
\]
the kernel \(K_n=\ker(M_{n-1}\to M_{n-2})\in\mathcal{G}_C^j.\)
\end{itemize}
Moreover, if \(\G_C^j\)-dim\(_R(M)<\infty,\) then \(\G_C^j\)-dim\(_R(M)=\alpha_R(M)-j.\)
\end{thm}

\pf: Clearly, the equivalence of \((i)\) and \((ii)\) will imply the last statement.
\begin{itemize}
\item[\((i)\Rightarrow (ii)\)] Suppose that \(\G_C^j\)-dim\(_R(M)\leq n\) and consider a \(\G_C^j\)-resolution of \(M\)
\[
\begin{tikzcd}
0 \arrow[r] & M_n \arrow[r] & M_{n-1} \arrow[r] & \cdots \arrow[r] & M_1 \arrow[r] & M_0 \arrow[r] & M \arrow[r] & 0
\end{tikzcd}
\]
By Lemma \ref{GCExtShiftLem} we have \(0=\Ext_R^{\g_R(M_n)+i}(M_n,C)\cong\Ext_R^{\g_R(M_n)+n+i}(M,C)\) for \(i>0\) since \(M_n\in\mathcal{G}_C^j.\)  Then
\[
\g_R(M_n)+n\geq \alpha_R(M)\Rightarrow n\geq \alpha_R(M)-\g_R(M_n)
\]
\item[\((ii)\Rightarrow (i)\)]  Suppose that \(G_C^j\)-dim\(_R(M)=p.\)  If \(p\leq n\) we are done, and so we may assume that \(p>n.\)  Consider a \(G_C^j\)-resolution of \(M\)
\[
\begin{tikzcd}
0 \arrow[r] & M_p \arrow[r] & \cdots \arrow[r] & M_n \arrow[r] & M_{n-1} \arrow[r] & \cdots \arrow[r] & M_1 \arrow[r] & M_0 \arrow[r] & M \arrow[r] & 0
\end{tikzcd}
\]
and let \(K_n=\ker(M_{n-1}\to M_{n-2}).\)  Then by Lemma \ref{GCExtShiftLem} again we have
\[
\Ext_R^{\g_R(M_{n-2})+i}(K_n,C)\cong\Ext_R^{\g_R(M_{n-1})+n+i}(M,C).
\]
Since \(n\geq \alpha_R(M)-j\) and \(j=\g_R(M_{n-1}),\) we have \(\g_R(M_{n-1})+n+i\geq\alpha_R(M)+i.\)  Thus\\ \(\Ext_R^{\g_R(M_{n-1})+i}(K_n,R)=0\) for \(i>0.\)  
Therefore \(\g_R(K_{n})=j=\alpha_R(K_n).\)  So by Proposition \ref{GCFinitedimPerfect0} we see that \(K_n\in\mathcal{G}_C^j.\)  Therefore \(G_C^j\)-dim\(_R(M)\leq n.\)
\item[\((i)\Rightarrow (iii)\)] Suppose \(G_C^j\)-dim\(_R(M)\leq n.\)  Consider a \(G_C^j\)-resolution of \(M\) of length \(n\)
\[
\begin{tikzcd}
0 \arrow[r] & M_n \arrow[r] & M_{n-1} \arrow[r] & \cdots \arrow[r] & M_0 \arrow[r] & M \arrow[r] & 0
\end{tikzcd}
\]
and let \(K_n=\ker(M_{n-1}\to M_{n-2}).\)  Consider a specific \(G_C^j\)-resolution constructed by using projective \(R/(\bar{x})\)-modules, where \(\bar{x}\) is a regular \(C\)-sequence in \(\ann_R(M)\) of length \(j,\)
\[
\begin{tikzcd}
0 \arrow[r] & S_n \arrow[r] & P_{n-1} \arrow[r] & \cdots \arrow[r] & P_0 \arrow[r] & M \arrow[r] & 0
\end{tikzcd}
\]
It is then sufficient to prove that \(S_n\in\mathcal{G}_C^j\) if and only if \(K_n\in\mathcal{G}_C^j.\)  Since \(P_i\) is projective (as \(R/(\bar{x})\)-modules) we get mappings from one complex to the other.  Then the result follows by considering the mapping cone and using Proposition \ref{GCShortexactsequence} \((a)\) and \((c).\)
\item[\((iii)\Rightarrow (i)\)] Clear.
\end{itemize}\qed

\begin{cor}\label{GCjinequality}
If \(C\) is a semidualizing \(R\)-module and \(M\) is a finitely generated \(R\)-module, then
\[
\G_C^i\text{-dim}_R(M)\leq \G_C^j\text{-dim}_R(M)
\]
for any \(i\geq j.\)  Moreover, if \(\G_C^j\)-dim\(_R(M)<\infty\) then \(\G_C^i\)-dim\(_R(M)=\G_C^j\)-dim\(_R(M)+(j-i)\) for all \(i\geq j.\)
\end{cor}

So we get this chain of inequalities
\[
\G_C^{\g_R(M)}\text{-dim}_R(M)\leq \G_C^{\g_R(M)-1}\text{-dim}_R(M)\leq\cdots\leq \G_C^1\text{-dim}_R(M)\leq \G_C^0\text{-dim}_R(M)
\]
In fact we get the following
\begin{thm}\label{GCdimallalwaysfinite}
Let \(C\) be a semidualizing \(R\)-module and \(M\) a finitely generated \(R\)-module.  Then \(\G_C^i\)-dim\(_R(M)<\infty\) if and only if \(\G_C^{\g_R(M)}\)-dim\(_R(M)<\infty\) for any \(i\leq \g_R(M).\)
\end{thm}

\pf: Clearly if \(\G_C^i\)-dim\(_R(M)<\infty,\) then \(\G_C^{\g_R(M)}\)-dim\(_R(M)<\infty\) by Corollary \ref{GCjinequality}.  So suppose that \(\G_C^{\g_R(M)}\)-dim\(_R(M)<\infty.\)  
Then for a suitable \(G_C^{\g_R(M)}\)-resolution of \(M\) (constructed in \(R/(\bar{x})\) for \(\bar{x}\) a regular \(C\)-sequence in \(\ann_R(M)\)) we can convert it into a \(\G_C\)-resolution of \(M\) 
by considering it as a complex over \(R/(\bar{x}).\)  As finiteness of \(\G_C\)-dimension is preserved between \(R\) and \(R/(\bar{x})\) when \(\bar{x}\) is a regular sequence 
\cite[Proposition 1.5.3]{GDIMChristensen} we are done by Corollary \ref{GCjinequality}.\qed

\begin{cor}\label{GCjshiftofGC}
Let $C$ be a semidualizing $R$-module and $M$ a finitely generated \(R\)-module.  If \(\G_C\)-dim\(_R(M)<\infty,\) then
\[
\G_C\text{-dim}_R(M)=\G_C^j\text{-dim}_R(M)+j
\]
\end{cor}

Note that from now on, any assumption about the finiteness of \(G_C\)-dimension is the same as assuming that any or all of the \(G_C^i\)-dimensions are finite.  
So in the rest of the results we will use notation that fits the theme of the result.  The next result is an Auslander-Bridger formula for these dimensions.

\begin{cor}[Auslander-Bridger Formula for \(\mathcal{G}_C^j\)]\label{GCjdimdepthequation}
Let \(C\) be a semidualizing module for a local ring \((R,\mathfrak{m})\) with \(\G_C^j\)-dim\(_R(M)<\infty\) for an \(R\)-module \(M.\)  Then
\[
\G_C^j\text{-dim}_R(M)=\depth(R)-\depth_R(M)+j
\]
\end{cor}

\section{Modified Serre Conditions and \(\mathbf{G_C^j}\)-dimension}

We will now recall notation and results from both \cite{LinkageModulesSerreSadDib} and \cite{QuasiPerfectModFoxby} to help present the next results.  The goal of this section is to provide results concerning Serre-like conditions using these newly developed dimensions.  These Serre-like conditions are then used to provide results concerning module linkage.\newline

Given a semidualizing \(R\)-module \(C\) and a finitely generated \(R\)-module \(M,\) one can take a projective presentation
\[
P_1 \to P_0 \to M \to 0
\]
and by applying \(\Hom_R(-,C)\) on obtains the exact sequence
\[
\begin{tikzcd}
0 \arrow[r] & \Hom_R(M,C) \arrow[r] & \Hom_R(P_0,C) \arrow[r] & \Hom_R(P_1,C) \arrow[r] & D_CM \arrow[r] & 0 
\end{tikzcd}
\]
where \(D_C(M)=\coker\!(\Hom_R(P_0,C)\to \Hom_R(P_1,C)).\)  In \cite{LinkageModulesSerreSadDib}, this is called the \(C\) transpose of \(M\) with notation \(\text{Tr}_CM.\)  Using the notation \((-)^{\nabla}=\Hom_R(-,C),\) one can write the following exact sequences which arise from the above sequence
\[
\begin{tikzcd}
0 \arrow[r] & \Ext_R^1(D_CM,C) \arrow[r] & M \arrow[r] & M^{\nabla\nabla} \arrow[r] & \Ext_R^2(D_CM,C) \arrow[r] & 0
\end{tikzcd}
\]
\[
\begin{tikzcd}
0 \arrow[r] & \Ext_R^1(M,C) \arrow[r] & D_CM \arrow[r] & (D_CM)^{\nabla\nabla} \arrow[r] & \Ext_R^2(M,C) \arrow[r] & 0
\end{tikzcd}
\]

Now suppose that \(\g_R(M)\geq j.\)  In the same manner we can take a \(\G_C^j\)-presentation of \(M,\) i.e. a presentation
\[
M_1 \to M_0 \to M \to 0
\]
where \(M_1,M_0\in\mathcal{G}_C^j\) and after applying \(\Ext_R^j(-,C)\) we get the exact sequence
\[
\begin{tikzcd}
0 \arrow[r] & \Ext_R^j(M,C) \arrow[r] & \Ext_R^j(M_0,C) \arrow[r] & \Ext_R^j(M_1,C) \arrow[r] & D_C^jM \arrow[r] & 0
\end{tikzcd}
\]
which leads to the following result.
\begin{lem}\label{GCddualexactsequence}
Let \(C\) be a semidualizing \(R\)-module and \(M\) an \(R\)-module with \(\g_R(M)\geq j.\)  Then there is an exact sequence
\[
\begin{tikzcd}
0 \arrow[r] & \Ext_R^{j+1}(D_C^gM,C) \arrow[r] & M \arrow[r, "\delta_C^j(M)"] & \Ext_R^j(\Ext_R^j(M,C),C) \arrow[r] & \Ext_R^{j+2}(D_C^gM,C) \arrow[r] & 0
\end{tikzcd}
\]
\end{lem}

\pf: Consider the presentation of \(M\) (as above)
\[
\begin{tikzcd}
(\Pi): & M_1 \arrow[r, "\alpha"] & M_0 \arrow[r, "\beta"] & M \arrow[r] & 0
\end{tikzcd}
\]
Applying \(\Ext_R^j(-,C)\) to this and splitting the resulting complex into short exact sequences we get
\[
\begin{tikzcd}
(\Pi_0): & 0 \arrow[r] & \Ext_R^j(M,C) \arrow[r, "\beta^*"] & \Ext_R^j(M_0,C) \arrow[r, "\gamma_0"] & Q \arrow[r] & 0\\
(\Pi_1): & 0 \arrow[r] & Q \arrow[r, "\gamma_1"] & \Ext_R^j(M_1,C) \arrow[r] & D_C^gM \arrow[r] & 0
\end{tikzcd}
\]
where \(\gamma_1\circ\gamma_0=\alpha^*.\)  Applying \(\Ext_R^j(-,C)\) to both of these gives the exact sequences
\[
\begin{tikzcd}[column sep = 2ex]
(\Pi_0^*): & 0 \arrow[r] & \Ext_R^j(Q,C) \arrow[r, "\gamma_0^*"] & \Ext_R^j(\Ext_R^j(M_0,C),C) \arrow[r, "\beta^{**}"] & \Ext_R^j(\Ext_R^j(M,C),C) \arrow[r] & \Ext_R^{j+1}(Q,C) \arrow[r] & 0\\
(\Pi_1^*): & 0 \arrow[r] & \Ext_R^j(D_C^gM,C) \arrow[r] & \Ext_R^j(\Ext_R^j(M_1,C),C) \arrow[r, "\gamma_1^*"] & \Ext_R^j(Q,C) \arrow[r] & \Ext_R^{j+1}(D_C^gM,C) \arrow[r] & 0.
\end{tikzcd}
\]
Then if we consider the following commutative diagram
\[
\begin{tikzcd}[row sep = 10ex]
            & M_1 \arrow[r, "\alpha"] \arrow[d, "\gamma_1^*\circ\delta_C^j(M_1)"] & M_0 \arrow[r, "\beta"] \arrow[d, "\delta_C^j(M_0)"] & M \arrow[r] \arrow[d, "\delta_C^j(M)"] & 0\\
0 \arrow[r] & \Ext_R^j(Q,C) \arrow[r, "\gamma_0^*"] & \Ext_R^j(\Ext_R^j(M_0,C),C) \arrow[r, "\beta^{**}"] & \Ext_R^j(\Ext_R^j(M,C),C).
\end{tikzcd}
\]
Now \(\coker\!(\gamma_1^*\circ\delta_C^j(M_1))=\coker\!(\gamma_1^*)=\Ext_R^{j+1}(D_C^gM,C)\) according to \((\Pi_1^*).\)  Therefore by the Snake Lemma, \(\ker(\delta_C^j(M))\cong \coker(\gamma_1^*\circ\delta_C^j(M_1))\cong \Ext_R^{j+1}(D_C^jM,C).\)\\

Further, \(\coker(\delta_C^g(M))\cong\coker(\beta^{**})\cong \Ext_R^{j+1}(Q,C)\cong\Ext_R^{j+2}(D_C^gM,C),\) where the last isomorphism comes from applying \(\Ext_R^j(-,C)\) to \((\Pi_1).\)\\

Lastly, we see from applying \(\Ext_R^j(-,C)\) to \((\Pi_0)\) and \((\Pi_1)\) that \(\Ext_R^{j+i}(\Ext_R^j(M,C),C)\cong \Ext_R^{j+1+i}(Q,C)\) and \(\Ext_R^{j+i}(Q,C)\cong \Ext_R^{j+i+1}(D_C^jM,C)\) for \(i>0,\) 
respectively.  This says that \(\Ext_R^{j+i}(\Ext_R^j(M,C),C)\cong \Ext_R^{j+i+2}(D_C^gM,C)\) for \(i>0,\) as desired.\qed

In the same way one obtains a sequence of functors
\[
\begin{tikzcd}
0 \arrow[r] & \Ext_R^{j+1}(D_C^gM,-) \arrow[r] & M\otimes_R- \arrow[r] & \Ext_R^j(\Ext_R^j(M,C),-) \arrow[r] & \Ext_R^{j+2}(D_C^gM,-) \arrow[r] & 0 
\end{tikzcd}
\]
for any semidualizing \(R\)-module \(C.\)

\begin{cor}\label{GCdimzeroifandonlyifAusdual}
Let \(C\) be a semidualizing \(R\)-module and \(M\) an \(R\)-module with \(\g_R(M)=g.\)  Then \(M\in\mathcal{G}_C^g\) if and only if \(\Ext_R^{g+i}(M,C)=\Ext_R^{g+i}(D_C^gM,C)=0\) for \(i>0.\)  
Consequently, \(M\in\mathcal{G}_C^g\) if and only if \(D_C^gM\in\mathcal{G}_C^g.\)
\end{cor}

\begin{cor}\label{GCAusduallongexact}
Suppose that \(0\to M'\to M\to M''\to 0\) is a short exact sequence with
\[
\min\{\g_R(M'),\g_R(M),\g_R(M'')\}\geq g.
\]
Then there is a long exact sequence
\[
0 \to \Ext_R^g(M'',C) \to \Ext_R^g(M,C) \to \Ext_R^g(M',C) \to D_C^gM'' \to D_C^gM \to D_C^gM' \to 0
\]
\end{cor}

With these results in hand we can begin to introduce Serre-like conditions.  Recall that an \(R\)-module \(M\) satisfies \((S_n)\) if \(\depth_{R_p}(M_p)\geq \min\{n,\dim(R_p)\}\) for \(p\in\text{Spec}(R).\)  
One says that an \(R\)-module \(N\) satisfies \(\tilde{S}_n\) if \(\depth_{R_p}(N_p)\geq \min\{n,\depth(R_p)\}\) for \(p\in\text{Spec}(R).\)  
It is natural for us to generalize \(\tilde{S}_n\) since we will want to use our Auslander-Bridger type formula which uses the depth of \(R\) and not the dimension.  
So, we define the following
\begin{defn}[Generalized Serre Condition]\label{GCjSerrecondition}
Let \(C\) be a semidualizing \(R\)-module and \(M\) an \(R\)-module.  We say that \(M\) satisfies \(\tilde{S}_n^g\) if
\[
\depth_{R_p}(M_p)+g\geq \min\{n,\depth(R_p)\}\quad\forall p\in\text{Spec}(R)
\]
\end{defn}

Note, \(\tilde{S}_n^g\) is always satisfied when \(n\leq g.\)  There is a similar generalized condition, \((S_n^j)\) for rings which has been shown to preserve many of the same results as for \((S_n)\) \cite{GenSerreHolmes}.  We also make the following definition

\begin{defn}\label{GCntorsionlessdefn}
Let \(C\) be a semidualizing \(R\)-module and \(M\) an \(R\)-module with \(\g_R(M)\geq g.\)  Then we say that \(M\) is \(\mathbf{C_n^g}\)\textbf{-torsionless} if \(\Ext_R^{g+i}(D_C^gM,C)=0\) for \(1\leq i\leq n.\) 
\end{defn}

With these definitions we prove Theorem \ref{GCSerreconditionthm} and Corollary \ref{GCSerreConditionCor}, generalizing \cite[Proposition 2.4]{LinkageModulesSerreSadDib} and \cite[Proposition 2.7]{LinkageModulesSerreSadDib}, respectively.  
The next proposition helps give intuition as to how these definitions fit together.  We show that a module being \(C_{n-g}^g\)-torsionless is stronger than satisfying \(\tilde{S}_n^g.\)

\begin{prop}\label{GCSerreconditionprop}
Let \(C\) be a semidualizing \(R\)-module and \(M\) an \(R\)-module with \(\g_R(M)=g\) and \(n\geq g,\) then for the following conditions:
\begin{itemize}
\item[\((i)\)] \(M\) is \(C_{n-g}^g\)-torsionless
\item[\((ii)\)] \(M\) is an \((n-g)^{\text{th}}\) \(\mathcal{G}_C^g\)-syzygy
\item[\((iii)\)] \(M\) satisfies \(\tilde{S}_n^g\)
\end{itemize}
we have \((i)\Rightarrow (ii) \Rightarrow (iii).\)
\end{prop}

\pf: In the following, we assume that \(n>g\) as the case $n=g$ is vacuously satisfied.
\begin{itemize}
\item[\((i)\Rightarrow(ii)\)] Consider a \(\mathcal{G}_C^g\)-resolution of \(\Ext_R^{g}(M,C)\)

\[
\begin{tikzcd}
\cdots \arrow[r] & M_{n-g-1} \arrow[r] & \cdots \arrow[r] & M_0 \arrow[r] & \Ext_R^{g}(M,C) \arrow[r] & 0
\end{tikzcd}
\]
If we apply \(\Ext_R^{g}(-,C)\) we get a complex
\[
\begin{tikzcd}
(*): & 0 \arrow[r] & \Ext_R^{g}(\Ext_R^{g}(M,C),C) \arrow[r] & \Ext_R^{g}(M_0,C) \arrow[r] & \cdots \arrow[r] & \Ext_R^{g}(M_{n-g-1},C)
\end{tikzcd}
\]
When \(n-g=1\) we have that \(M\subset \Ext_R^{g}(\Ext_R^{g}(M,C),C)\) by Lemma \ref{GCddualexactsequence} and are done using \((*).\)  
If \(n-g=2,\) then \(M\cong \Ext_R^{g}(\Ext_R^{g}(M,C),C)\) by Lemma \ref{GCddualexactsequence} and it is clear that \(M\) is a \(2^{\text{nd}}\) \(\mathcal{G}_C^g\)-syzygy.  
So suppose that \(n-g>2.\)  Then as \(\Ext_R^{g+i}(D_C^gM,C)\cong \Ext_R^{g+i-2}(\Ext_R^{g}(M,C),C)\) for \(2<i\leq n-g\) we see that \((*)\) is exact.

\item[\((ii)\Rightarrow (iii)\)] If we take an exact complex
\[
\begin{tikzcd}
0 \arrow[r] & M \arrow[r] & M_1 \arrow[r] & \cdots \arrow[r] & M_{n-g}
\end{tikzcd}
\]
then clearly \(\depth_{R_p}(M_p)\geq \min\{\depth_{R_p}\left((M_{n-g})_p\right),n-g\}.\)  Using Corollary \ref{GCjdimdepthequation} we have
\[
\depth_{R_p}(M_{n-g})=\depth(R_p)-\g_{C_p}(M_p)
\]
and so
\[
\depth_{R_p}(M_p)\geq \min\{\depth(R_p)-\g_{C_p}(M_p),n-g\}\geq\min\{\depth(R_p),n\}-\g_{C_p}(M_p)
\]
\end{itemize}\qed

We can extend this to an equivalence of statements assuming that the \(\G_C^g\)-dimension of \(M\) is locally finite.
For notational convenience we adopt the following notation
\[X^{i+n-1}=\{p\in\text{Spec}(R)~|~ \depth(R_p)\leq i+n-1\}.\]

\begin{thm}\label{GCSerreconditionthm}
Let \(C\) be a semidualizing \(R\)-module, \(M\) an \(R\)-module with \(\g_R(M)=g,\) \(n\geq g,\) and \(M\) having locally finite $G^g_C-$dimension. Then the following are equivalent:
\begin{itemize}
\item[\((i)\)] \(M\) is \(C_{n-g}^g\)-torsionless
\item[\((ii)\)] \(M\) is an \((n-g)^{\text{th}}\) \(\mathcal{G}_C^g\)-syzygy
\item[\((iii)\)] \(M\) satisfies \(\tilde{S}_n^g\)
\item[\((iv)\)] \(\g_{C_p}(\Ext_{R_p}^{\g_{C_p}(M_p)+i}(M_p,C_p))\geq i+n\) for \(i\geq 1\) and \(p\in X^{i+n-1}\)
\end{itemize}
\end{thm}

\pf: We have already seen that \((i)\Rightarrow (ii) \Rightarrow (iii).\)  We will now show that \((iii)\Rightarrow (iv)\) and \((iv)\Rightarrow (i).\)
\begin{itemize}
\item[\((iii)\Rightarrow (iv)\)] Fix \(i\geq 1\) and a prime \(p\in\text{Spec}(R)\) with \(\depth(R_p)<i+n.\)  Let \(\g_{C_p}(M_p)=g_p.\)  We need to show that \(\Ext_{R_p}^{g_p+i}(M_p,C_p)=0.\)  
To do this we will show that \(p\notin \Supp(\Ext_{R_p}^{g_p+i}(M_p,C_p)).\)  Using Corollary \ref{GCjdimdepthequation} we have
\begin{align*}
\G_{C_p}^{g_p}\text{-dim}_{R_p}(M_p)&=\depth(R_p)-\depth_{R_p}(M_p)-g_p\\
&\leq \depth(R_p)-\min\{n, \depth(R_p)\}\\
&=\max\{0,\depth(R_p)-n\}
\end{align*}
and since \(\depth(R_p)<i+n\) we have \(\G_{C_p}^{g_p}\)-dim\(_{R_p}(M_p)<i+n-n=i.\)  This says that \(\alpha_{C_p}(M_p)<i+g_p\) and so \(\Ext_{R_p}^{g_p+i}(M_p,C_p)=0\)
\item[\((iv)\Rightarrow (i)\)]  It is enough to show this in a local ring \((R,p)\) for \(p\in\text{Spec}(R).\)  So we may assume, \(\G_C^g\)-dim\(_R(M)=\alpha_R(M)-g<\infty.\)  
By Corollary \ref{GCdimzeroifandonlyifAusdual} the result holds if \(\G_C^g$-dim$_R(M)=0.\)  So suppose that \(\G_C^g\)-dim\(_R(M)=p>0\) and we proceed by induction on \(p.\)  
Let
\[
0 \to K \to N \to M \to 0
\]
be a short exact sequence with \(N\in\mathcal{G}_C^g\) and so \(\G_C^g\)-dim\(_R(K)=p-1.\)  Then we have \(\g_R(\Ext_R^{g+i}(K,C))=\g_R(\Ext_R^{g+i+1}(M,C))\geq i+n+1\) for \(i>0.\)  
So by induction \(\Ext_R^{g+i}(D_C^gK,C)=0\) for \(1\leq i\leq n-g+1.\)  That is \(K\) is \(C_{n-g+1}^g\)-torsionless.  The result will now follow by the following Lemma.
\end{itemize}\qed

\begin{lem}\label{}
Let \(0 \to M' \to M \to M'' \to 0\) be a short exact sequence with
\[
\min\{\g_R(M'),\g_R(M),\g_R(M'')\}\geq g
\]
and \(Q=\coker\!(\Ext_R^g(M,C)\to\Ext_R^g(M',C)).\)  If \(M'\) is \(C_{k+1}^g\)-torsionless, \(M\) is \(C_k^g\)-torsionless, and \(\g_R(Q)\geq k+1\) then \(M''\) is \(C_k^g\)-torsionless.
\end{lem}

\pf: From Corollary \ref{GCAusduallongexact} we have the exact sequence
\[
\begin{tikzcd}
0 \arrow[r] & Q \arrow[r] & D_C^gM'' \arrow[r] & D_C^gM \arrow[r] & D_C^gM' \arrow[r] & 0.
\end{tikzcd}
\]
Breaking this into short exacts sequences and looking at the corresponding long exact sequences in Ext gives the desired result.\qed

The following corollary will be useful for proofs in the next sections.

\begin{cor}\label{GCSerreConditionCor}
Let \(C\) be a semidualizing \(R\)-module with locally finite injective dimension, \(M\) an \(R\)-module with \(\g_R(M)=g,\) \(n\geq g,\)   Then the following are equivalent:
\begin{itemize}
\item[\((i)\)] \(M\) is \(C_{n-g}^g\)-torsionless
\item[\((ii)\)] \(M\otimes_R C\) is an \((n-g)^{\text{th}}\) \(\mathcal{G}_C^g\)-syzygy
\item[\((iii)\)] \(M\otimes_R C\) satisfies \(\tilde{S}_n^g\)
\item[\((iv)\)] \(\g_{C_p}(\Ext_{R_p}^{\g_{R_p}(M_p)+i}(M_p,R_p))\geq i+n\) for \(i> 0\) and \(p\in X^{i+n-1}(R).\)
\end{itemize}
\end{cor}

\pf: Let \(N=M\otimes_R C.\)  Since \(C\) is semidualizing
\[
\Ext_R^g(\Ext_R^g(N,C),C) = \Ext_R^g(\Ext_R^g(M\otimes_R C,C),C)\cong \Ext_R^g(\Ext_R^g(M,R),C).
\]
Then using the exact sequence from Lemma \ref{GCddualexactsequence} and the remark following the lemma we have the following commutative diagram
\[
\begin{tikzcd}
0 \arrow[r] & \Ext_R^{g+1}(D_R^gM,C) \arrow[r] \arrow[d] & N \arrow[r] \arrow[d, equals] & \Ext_R^g(\Ext_R^g(M,R),C) \arrow[r] \arrow[d, "\cong"] & \Ext_R^{g+2}(D_R^gM,C) \arrow[r] \arrow[d] & 0\\
0 \arrow[r] & \Ext_R^{g+1}(D_C^gN,C) \arrow[r] & N \arrow[r] & \Ext_R^g(\Ext_R^g(N,C),C) \arrow[r] & \Ext_R^{g+2}(D_C^gN,C) \arrow[r] & 0
\end{tikzcd}
\]
So it follows naturally that \(\Ext_R^{g+i}(D_R^gM,C)\cong \Ext_R^{g+i}(D_C^gN,C)\) for \(i>0\) using the above diagram and \(\Ext_R^g(M,R)\cong\Ext_R^g(N,C).\)  
Then the result follows by replacing \(M\) with \(N\) in Theorem \ref{GCSerreconditionthm}.\qed

To close this section, we will prove a result that is analogous to \cite[Theorem 2.12]{LinkageModulesSerreSadDib}.  
We first need to recall the concept of the Auslander class of a semidualizing module.  
This module class were defined by Foxby \cite{QuasiPerfectModFoxby} and further developed by Avramov and Foxby in \cite{RingHomoGDimAvramovFoxby}.

\begin{defn}\label{AuslanderClassDefn}
Let \(C\) be a semidualizing \(R\)-module.  The \textbf{Auslander class with respect to \(\mathbf C\)}, \(\mathcal{A}_C,\) consists of all \(R\)-modules \(M\) satisfying:
\begin{itemize}
\item[\((i)\)] The map \(M\to\Hom_R(C,M\otimes_R C)\) is an isomorphism
\item[\((ii)\)] \(\Tor_i^R(M,C)=0=\Ext_R^i(C,M\otimes_RC)\) for all \(i>0.\)
\end{itemize}
\end{defn}

Then using \cite[Lemma 2.11]{LinkageModulesSerreSadDib} we get

\begin{prop}\label{GCSerreConditionAll}
Let \(C\) be a semidualizing \(R\)-module, \(M\in\mathcal{A}_C\), \(\g_R(M)=g,\) \(n\geq g,\) and M has localy finite $G^g_C-$dimension.  Then the following are equivalent:
\begin{itemize}
\item[\((i)\)] \(M\) is \(R_{n-g}^g\)-torsionless
\item[\((ii)\)] \(M\) is \(C_{n-g}^g\)-torsionless
\item[\((iii)\)] \(M\otimes_R C\) satisfies \(\tilde{S}_n^g\)
\item[\((iv)\)] \(M\) satisfies \(\tilde{S}_n^g.\)
\end{itemize}
\end{prop}

\pf: The equivalence of \((iii)\) and \((iv)\) follows from \cite[Lemma 2.11]{LinkageModulesSerreSadDib} and that \(\g_R(M)=\g_R(M\otimes_R C)\) since
\[
\Ext_R^i(M\otimes_R C,C)\cong\Ext_R^i(M,\Hom_R(C,C))\cong \Ext_R^i(M,R)
\]
and \(\g_R(M)=\g_C(M).\)  Note that \((i)\) and \((iv)\) are equivalent from Theorem \ref{GCSerreconditionthm} by replacing \(C\) with \(R.\)  Further note that \((ii)\) and \((iii)\) are equivalent from Corollary \ref{GCSerreConditionCor}.\qed


\section{Module Linkage}

Throughtout this section \(R\) will be a semiperfect ring.  Semiperfect rings are products of commutative local Noetherian rings, and so we do not lose too much generality (with regards to the rest of the paper) in this setting.  One of the reasons for choosing such a setting is that in a semiperfect ring every finitely generated module has a projective cover, (see \cite{AlgebraIIFaith} Ch. 18).  In fact, we will only need this hypothesis for Proposition \ref{horizontallylinkedimpliesgstable} and the results following it.\newline

The goal of this section is to use intermediate \(C\)-Gorenstein dimensions to help extend some results about module linkage that are found in \cite{LinkageMartsinkovskyStrooker, LinkageModulesSerreSadDib}.  To repeat what has been said before, given an \(R\)-module \(M\) of grade zero we can take a projective presentation
\[
\begin{tikzcd}
P_1 \arrow[r] & P_0 \arrow[r] & M \arrow[r] & 0
\end{tikzcd}
\]
and get the following exact sequence
\[
\begin{tikzcd}
0 \arrow[r] & \Hom_R(M,R) \arrow[r] & P_0^* \arrow[r] & P_1^* \arrow[r] & D_RM \arrow[r] & 0
\end{tikzcd}
\]
We then set 
\[
\lambda_RM=\Omega(D_RM)=\ker(P_1^*\to D_RM)=\coker(\Hom_R(M,R)\to P_0^*).
\]
A module \(N\) is said to be \textit{directly linked to \(M\)} if \(M\cong \lambda_RN\) and \(N\cong\lambda_RM.\)  A module \(M\) is \textit{horizontaly linked} (to \(\lambda_RM\)) if \(M\cong\lambda_R^2M.\)    
One can then prove that a module is horizontally linked if and only if it is stable and \(R_1^0\)-torsionless, \cite[Theorem 2]{LinkageMartsinkovskyStrooker}.  
There are many nice results about this linkage and one that we single out is that this linkage preserves \(C\)-Gorenstein dimension zero.  
This is completely expected as this linkage was constructed to emulate complete intersection ideal linkage. 
We also reiterate the more general notion of module linkage, due to Nagel (see \cite{NagelLiaison}).\newline

Given a semidualizing \(R\)-module \(C,\) an \(R\)-module \(Q\) is said to be \textbf{C-quasi-Gorenstein} if \(Q\in\mathcal{G}_C^{\g_R(Q)}\) and there is some isomorphism \(\alpha : Q\to \Ext_R^{\g_R(Q)}(Q,C).\)  
\(C$-Quasi-Gorenstein modules, \(Q,\) have nice properties \cite{SturgeonQGorenstein,BrennanYorkQGorenstein,BrennanYorkETFS,NagelLiaison}.  
Given a \(C\)-quasi-Gorenstein module \(Q\) of grade \(g\) we will denote by \(\Epi(Q)\) the set of all $R$-module homomorphisms \(\varphi:Q\to M\) where \(\im\varphi\) has the same grade as \(Q.\)  
Given such a homomorphism \(\varphi\) we have a short exact sequence
\[
\begin{tikzcd}
0 \arrow[r] & \ker\varphi \arrow[r] & Q \arrow[r] & \im\varphi \arrow[r] & 0
\end{tikzcd}
\]
which induces a long exact sequence
\[
\begin{tikzcd}
0 \arrow[r] & \Ext_R^g(\im\varphi,C) \arrow[r] & \Ext_R^g(Q,C) \arrow[r, "\psi"] & \Ext_R^g(\ker\varphi,C) \arrow[r] & \Ext_R^{g+1}(\im\varphi,C) \arrow[r] & \cdots
\end{tikzcd}
\]
If \(\alpha:Q\to\Ext_R^g(Q,C)\) is an isomorphism, then we can construct a short exact sequence from the long sequence above as
\[
\begin{tikzcd}
0 \arrow[r] & \Ext_R^g(\im\varphi,C) \arrow[r] & Q \arrow[r, "L_Q(\varphi)"] & \im L_Q(\varphi) \arrow[r] & 0
\end{tikzcd}
\]
where \(L_Q(\varphi)=\psi\circ\alpha.\)  

\begin{defn}\label{ModuleLinkageDefn}
We say that \(M\) and \(N\) are \textbf{directly linked by the } \(\mathbf{C}\)\textbf{-quasi-Gorenstein module } \(\mathbf{Q}\) if there are \(\varphi,\psi\in\Epi(Q)\) such that
\begin{enumerate}
\item[\((i)\)] \(M=\im\varphi,\) \(N=\im\psi\)
\item[\((ii)\)] \(M\cong \im L_Q(\psi),\) \(N\cong \im L_Q(\varphi).\)
\end{enumerate}
\end{defn}
It may not immediately be apparent that this is a generalization of the module linkage of Martsinkovsky and Strooker. If one restricts the quasi-Gorenstein modules by which one is allowed to be link itbecomes clear, see \cite[Remark 3.20]{NagelLiaison}.  We will also say that \(M\) and \(N\) are in the same linkage class if they are linked together by a chain of direct links.\newline  

With these definitions and results presented above it is easy to ask when we would run across rings with semidualizing modules \(C\) such that there are nontrivial \(C\)-quasi-Gorenstein modules.  Then the results that follow would be useful as they would give information about such modules that was not previously apparent.\newline

The following construction which is an extension of the ideas of \cite{SatherWagstaffNotes} allows one to build a family of such rings.  

\begin{exam}\label{semidualizingexampleconstruction}
Let \(k\) be a field and consider the ring \(R=k[X,Y]/(X,Y)^2.\)  Then $R$ is a local Cohen-Macaulay ring which is not Gorenstein.  It is clear then that \(R\) is a free \(k\)-module.  \(R\) has two semidualizing modules, 
\(R\) and \(\omega_R=\Hom_k(R,k).\)  In fact, \(\omega_R\) is dualizing and \(\omega_R\not\cong R.\)\newline

Clearly, \(R\) is an \(R\)-quasi-Gorenstein module.  It is then straightforward to see that \(k\) is a \(\omega_R\)-quasi-Gorenstein module as
\[
\Hom_R(k,\omega_R)=\Hom_R(k,\Hom_k(R,k))\cong \Hom_k(k\otimes_R R,k)\cong \Hom_k(k,k)\cong k
\]
Therefore, one has that \(R\) is an \(\omega_R\)-quasi-Gorenstein module, as \(R\) is a free \(k\)-module of rank 3.  
Therefore \(k\) and \(R\) are directly linked by the \(\omega_R\)-quasi-Gorenstein module \(k\oplus R.\)  However, \(k\) is not an \(R\)-quasi-Gorenstein module and so it is not easy to see if \(k\) and \(R\) are linked 
using some \(R\)-quasi-Gorenstein module.  Thus, \(k\) and \(R\) may not be directly linked through some \(R\)-quasi-Gorenstein module, but they are through some \(\omega_R\)-quasi-Gorenstein module.  
So we can see that the extension of linkage to semidualizing modules has given new information about linkage classes.\newline

Now, if we take \((R,m,k)\) as above and construct \(S=R[U,V]/(U,V)^2\) (just as \(R\) is constructed above) then \(S\) is a local Cohen-Macaulay ring with residue field \(k.\)  Then \(S\) has four distinct semidualizing modules \[S, \quad C_1=\Hom_R(S,R),\quad C_2=S\otimes_R \omega_R, \quad \hbox{\rm and} \omega_S=\Hom_R(S,\omega_R).\]  In fact, \(\omega_S\) is dualizing and \(\omega_S\not\cong S.\)  Once again, using \(S\)-quasi-Gorenstein modules it is not clear if \(S,\) \(R,\) and \(k\) are directly linked to each other.  However, \(R\) is a \(C_1\)-quasi-Gorenstein module as
\[
\Hom_S(R,C_1)=\Hom_S(R,\Hom_R(S,R))\cong\Hom_R(R\otimes_S S, R)\cong\Hom_R(R,R)\cong R
\]
and \(k\) is a \(\omega_S\)-quasi-Gorenstein module as
\[
\Hom_S(k,\omega_S)=\Hom_S(k,\Hom_R(S,\omega_R))\cong\Hom_R(k\otimes_S S,\omega_R)\cong\Hom_R(k,\omega_R)\cong k.
\]
Then, as \(S\) is a free \(R\)-module (of rank 3) we have that \(S\) and \(R\) are directly linked by the \(C_1\)-quasi-Gorenstein module \(S\oplus R.\)  Further, \(S\) is a free \(k\)-module (of rank 9) and so \(S\) and \(k\) are directly linked by the $\omega_S$-quasi-Gorenstein module $S\oplus k$.  As above, $R$ and $k$ are directly linked as well.
\end{exam}

One could continue such a construction and obtain a ring with \(2^n\) distinct semidualizing modules for any \(n.\)  
The idea to take away from this example is that by using different semidualizing modules for linkage  one to shift between the linkage classes of modules.  
So in what follows we are not only presenting new results concering linkage of modules but also presenting it in a way that concerns all such linkage classes with different semidualizing modules.\newline

We will now use the notation \(\mathcal{L}_Q(M)=\im L_Q(\varphi)\) where \(\varphi\in\Epi(Q)\) and \(\im\varphi=M.\)  Note that \(M\) is not necessarily directly linked to \(\mathcal{L}_Q(M).\)  
Nevertheless, we will call \(\mathcal{L}_Q(M)\) the \textit{first link} of \(M\) by \(Q.\)  The following will begin to illuminate some interactions between linkage and \(G_C^g\)-dimension.

\begin{prop}\label{pdimzeropreservedunderlinkage}
Suppose that \(M\in\mathcal{G}_C^g\) and \(M=\im\varphi\) for \(\varphi\in\Epi(Q).\)  Then \(\mathcal{L}_Q(M)\in\mathcal{G}_C^g.\) 
\end{prop}

\pf: We have the short exact sequence
\[
\begin{tikzcd}
0 \arrow[r] & \ker\varphi \arrow[r] & Q \arrow[r, "\varphi"] & M \arrow[r] & 0
\end{tikzcd}
\]
which induces the short exact sequence
\[
\begin{tikzcd}
0 \arrow[r] & \Ext_R^g(M,C) \arrow[r] & \Ext_R^g(Q,C) \arrow[r] & \Ext_R^g(\ker\varphi,C) \arrow[r] & 0.
\end{tikzcd}
\]
From the first sequence we see that \(\ker\varphi\in\mathcal{G}_C^g\) as both \(M\) and \(Q\) are in \(\mathcal{G}_C^g\) by Proposition \ref{GCShortexactsequence} \((a).\)  
Therefore \(\Ext_R^g(\ker\varphi,C)\in\mathcal{G}_C^g.\)  From the second sequence we see that \(\mathcal{L}_Q(M)\cong\Ext_R^g(\ker\varphi,C).\)\qed

\begin{lem}\label{pdimlinkagequotient}
Suppose that \(M\) is a \(C\)-quasi-Gorenstein \(R\)-module, \(N\subset M\) with \(N\in\mathcal{G}_C^{\g_R(M)},\) and \(\g_R(M/N)=\g_R(M).\)  
Then \(M/N\in\mathcal{G}_C^{\g_R(M/N)}\) if and only if \(M/N\) and \(\Ext_R^{\g_R(M)}(N,C)\) are directly linked by \(M.\)
\end{lem}

\pf: Let \(\g_R(M)=m.\)  Then the result is clear by considering the two exact sequences
\[
\begin{tikzcd}
0 \arrow[r] & N \arrow[r] & M \arrow[r] & M/N \arrow[r] & 0\
\end{tikzcd}
\]
\[
\begin{tikzcd}
0 \arrow[r] & \Ext_R^m(M/N,C) \arrow[r] & \Ext_R^m(M,C) \arrow[r] & \Ext_R^m(N,C) \arrow[r] & \Ext_R^{m+1}(M/N,C) \arrow[r] & 0
\end{tikzcd}
\]\qed

What is happening here is we are emulating the same idea as other linkage theories, but by considering a different ``presentation" of $M$ we get different information about the links.  
Since we are using a linking module of the same grade as our module, we can expect the linkage to preserve these types of properties.  In fact it is easily seen that the grade is preserved.  
We will follow the ideas in \cite{LinkageMartsinkovskyStrooker} to develop some properties of modules in \(\mathcal{G}_C^g\) under linkage.

\begin{defn}\label{gstabledefn}
We will say that \(M\) is \textbf{\(\mathcal{G}_C^g\)-stable} if it has no summands in \(\mathcal{G}_C^g.\)
\end{defn}

\begin{defn}\label{horizontallylinkedbyC}
Suppose that \(M\) is an \(R\)-module of grade \(g\) and \(\varphi\in\text{Epi}(Q)\) with \(\im\!(\varphi)=M.\)  
We say that \(M\) is \textbf{horizontally-linked by \(\mathbf{Q}\)} if \(M\cong\mathcal{L}^2_Q(M):=\mathcal{L}_Q(\mathcal{L}_Q(M)).\)
\end{defn}

It is implied in the definition that \(Q\) is a \(C\)-quasi-Gorenstein module.  
Now consider a \(\mathcal{G}_C^g\)-stable module \(M\) with \(M=\im\!(\varphi)\), \(\varphi\in\text{Epi}(Q)\) and \(\alpha:Q\to \Ext_R^g(Q,C)\) an isomorphism.  We can then dualize
\begin{center}
\begin{tikzcd}
0 \arrow[r] & \ker\varphi \arrow[r] & Q \arrow[r, "\varphi"] & M \arrow[r] & 0
\end{tikzcd}
\end{center}
and obtain
\begin{center}
\begin{tikzcd}
0 \arrow[r] & \Ext_R^g(M,C) \arrow[r, "\alpha^*\varphi^*"] & Q \arrow[r, "\psi\alpha"] & \mathcal{L}_Q(M) \arrow[r] & 0.
\end{tikzcd}
\end{center}
By dualizing once again we get
\begin{center}
\begin{tikzcd}
0 \arrow[r] & \Ext_R^g(\mathcal{L}_Q(M),C) \arrow[r, "\alpha^*\psi^*"] & Q \arrow[r, "\varphi^{**}\alpha^{**}"] & \mathcal{L}_Q^2(M) \arrow[r] & 0
\end{tikzcd}
\end{center}
Now it becomes clear that \(\im\!(\varphi^{**})\cong\mathcal{L}_Q^2(M)\) (\(\alpha^{**}\) is exactly \(\delta_C^g(Q^*)\alpha(\delta_C^g(Q))^{-1}\)).  
We may interchange \((-)^*\) and \(\Ext_R^g(-,C)\) when it is clear what grade we are in.  From the commutative square
\begin{center}
\begin{tikzcd}
Q \arrow[r, "\varphi", two heads] \arrow[d, "\delta_C^g(Q)"] & M \arrow[d, "\delta_C^g(M)"]\\
Q^{**} \arrow[r, "\varphi^{**}"] & M^{**}
\end{tikzcd}
\end{center}
we see that \(\im\!(\varphi^{**})\cong\im\!(\delta_C^g(M))\) when \(M\) is \(\mathcal{G}_C^g\)-stable.  This is summarized in the following result

\begin{prop}\label{imagesinlinkageandstuff}
Suppose that \(M\) is a \(\mathcal{G}_C^g\)-stable \(R\)-module with \(M=\im\!(\varphi)\), \(\varphi\in\text{Epi}(Q).\)  Then \(\im\!(\delta_C^g(M))\cong\mathcal{L}_Q^2(M)\) and we have a short exact sequence
\begin{center}
\begin{tikzcd}
0 \arrow[r] & \Ext_R^{g+1}(D_C^gM, C) \arrow[r] & M \arrow[r, "\delta_C^g(M)"] & \mathcal{L}_Q^2(M)\arrow[r] & 0.
\end{tikzcd}
\end{center}
\end{prop}
\qed
We see that for a \(\mathcal{G}_C^g\)-stable \(R\)-module \(M,\) \(\mathcal{L}_Q^2(M)\) is independent of \(Q.\)  In this case we drop the \(Q\) and write \(\mathcal{L}^2(M)\) and say that \(M\) is \textbf{horizontally-linked} if \(M\cong\mathcal{L}^2(M).\)  Horizontally linked modules are nice in the following way

\begin{prop}\label{horizontallylinkedimpliesgstable}
Suppose that \(M\) is horizontally-linked by \(Q\) and \(\g_R(M)=g.\)  Then \(M\) is \(\mathcal{G}_C^g\)-stable.
\end{prop}

\pf:  Suppose that \(M\cong\mathcal{L}^2_Q(M)\) and \(M\cong M'\oplus P\) where \(M'\) is \(\mathcal{G}_C^g\)-stable and \(P\in\mathcal{G}_C^g.\)  
As \(Q\to M\) is the start of a minimal \(\mathcal{G}_C^g\)-presentation, we get that $Q^*\to \mathcal{L}_Q(M)$ is one as well since $M'$ is $\mathcal{G}_C^g$-stable.  
Then we have
\[
Q\oplus P\twoheadrightarrow Q\cong Q^{**}\to\mathcal{L}_Q^2(M)\cong M\cong M'\oplus P
\]
Which gives a surjective endomorphism \(Q\oplus P\to Q\) (as \(\mathcal{L}_Q^2(M)\cong M\)) which must be an isomorphism.  Thus \(P=0.\)
\qed

\begin{prop}\label{linkishorizontallylinked}
Suppose that \(M\) is horizontally linked by \(Q.\)  Then so is \(\mathcal{L}_Q(M)\) and in particular \(\mathcal{L}_Q(M)\) is \(\mathcal{G}_C^{\g_R(M)}\)-stable.
\end{prop}

\pf: Since \(M\cong \mathcal{L}^2_Q(M)\) we know that \(\mathcal{L}^2_Q(\mathcal{L}_Q(M))\cong \mathcal{L}_Q(\mathcal{L}_Q^2(M))\cong\mathcal{L}_Q(M).\)\qed

These results lead to the characterization of horizontally-linked modules
\begin{thm}\label{horizontallylinkedtheorem}
A finitely generated \(R\)-module \(M\) of grade \(g\) is horizontally-linked if and only if it is \(\mathcal{G}_C^g\)-stable and \(C_{g+1}^g\)-torsionless.
\end{thm}

\pf: If \(M\) is horizontally-linked, then by Proposition \ref{horizontallylinkedimpliesgstable} \(M\) is \(\mathcal{G}_C^g\)-stable and by Proposition \ref{imagesinlinkageandstuff} \(M\) is \(C_{g+1}^g\)-torsionless.  The converse is clear by Proposition \ref{imagesinlinkageandstuff}.\qed

With this result in hand, we now turn to understanding how module linkage and \(G_C^g\)-dimension are related.  
We will prove Proposition \ref{LinkageandSerreLocalDuality1} and Corollary \ref{LinkageandSerreLocalDuality2} which are generalizations of Proposition 2.6 and Corollary 2.8 in \cite{LinkageModulesSerreSadDib}, respectively.  
First, recall the Local Duality Theorem \cite[Corollary 3.5.9]{BH}

\begin{thm}\label{LocalDualityTheorem}
Let \((R,\mathfrak{m},k)\) be a Cohen-Macaulay local ring of dimension \(d\) with a canonical module \(\omega_R.\)  Then for all finitely generated \(R\)-modules \(M\) and all integers \(i\) there exists natural isomorphism
\[
H_{\mathfrak{m}}^i(M)\cong\Hom_R(\Ext_R^{d-i}(M,\omega_R),E_R(k)),
\]
where \(E_R(k)\) is the injective envelope of \(k.\)
\end{thm}

At this point we relate the local cohomology of \(M\otimes_R \omega_R\) in a Cohen-Macaulay local ring to properties of the first link of \(M.\)

\begin{prop}\label{LinkageandSerreLocalDuality1}
Let \(R\) be a Cohen-Macaulay local ring of dimension \(d\) with canonical module \(\omega_R\) and \(M\) an \(R\)-module.  
Suppose that \(M\) is horizontally linked by \(Q\) with \(\g_R(M)=g,\) and that \(\Ext_R^{g+1}(D_R^gM,\omega_R)=0.\)  Then, for a positive integer \(n\geq g,\) the following statements are equivalent:
\begin{itemize}
\item[\((i)\)] \(\mathcal{L}_Q(M)\) satisfies \(\tilde{S}_n^g\)
\item[\((ii)\)] \(H_{\mathfrak{m}}^i(M\otimes_R\omega_R)=0\) for all \(i, d-n+g<i<d.\)
\end{itemize}
\end{prop}

\pf: \(\mathcal{L}_Q(M)\) satisfies \(\tilde{S}_n^g\) if and only if \(D_R^gM\) satisfies \(\tilde{S}_{n-1}^g\) because \(\mathcal{L}_Q(M)\) is a first \(\mathcal{G}_R^g\)-syzygy of \(D_R^gM\) and \(\Ext_R^{g+1}(D_R^gM,\omega_R)=0.\)  Thus by Theorem \ref{GCSerreconditionthm} statement \((i)\) is equivalent to \(D_R^gM\) being \((\omega_R)_{n-g-1}^g\)-torsionless, i.e.
\[
\Ext_R^{g+i}(D_{\omega_R}^g(D_R^gM),\omega_R)=0\text{ for all }i, 1\leq i\leq n-g-1.
\]
However, \(D_{\omega_R}^g(D_R^gM)\cong D_R^gD_R^gM\otimes_R\omega_R\) and as \(M\) is \(\mathcal{G}_R^g\)-stable, \(D_R^gD_R^gM\cong M\) and so \(D_R^gD_R^gM\otimes_R\omega_R\cong M\otimes_R\omega_R.\)  Hence \(\mathcal{L}_Q(M)\) satisfies \(\tilde{S}_n^g\) if and only if \(\Ext_R^{g+i}(M\otimes_R\omega_R,\omega_R)=0\) for all \(i, 1\leq i\leq n-g-1,\) which is equivalent to
\[
H_\mathfrak{m}^i(M\otimes_R\omega_R)=0\text{ for all }i, d-n+g<i<d
\]
by the Local Duality Theorem.\qed

\begin{cor}\label{LinkageandSerreLocalDuality2}
Let \(R\) be a Cohen-Macaulay local ring of dimension \(d\) with canonical module \(\omega_R\) and \(M\) an \(R\)-module.  Suppose that \(M\) is horizontally linked by \(Q\) with \(\g(M)=g,\) and that \(\Ext_R^{g+1}(D_R^gM,\omega_R)=0.\)  Then, for a positive integer \(n\geq g,\) the following statements are equivalent:
\begin{itemize}
\item[\((i)\)] \(M\otimes_R\omega_R\) satisfies \(\tilde{S}_n^g\)
\item[\((ii)\)] \(H_{\mathfrak{m}}^i(\mathcal{L}_Q(M))=0\) for \(d-n+g<i<d.\)
\end{itemize}
\end{cor}

\pf: This is clear using Proposition \ref{LinkageandSerreLocalDuality1} and Local Duality.\qed

Using Proposition \ref{GCSerreConditionAll} and Theorem \ref{horizontallylinkedtheorem} we get the following generalization of Corollary 2.14 in \cite{LinkageModulesSerreSadDib}.

\begin{thm}\label{GCDimandLinkageThm}
Let \(R\) be a Cohen-Macaulay ring, \(C\) a semidualizing \(R\)-module, and \(M\) a \(\mathcal{G}_C^g\)-stable \(R\)-module with \(\g_R(M)=g.\)  
Suppose that \(n\geq g,\) \(M\in\mathcal{A}_C,\) and the \(G_C^g-\hbox{dim}_R\) of  \(M)\) is locally finite. Then the following are equivalent:
\begin{itemize}
\item[\((i)\)] \(M\) satisfies \(\tilde{S}_n^g\)
\item[\((ii)\)] \(M\) is horizontally linked by some \(C\)-quasi-Gorenstein \(R\)-module \(Q\) and \(\Ext_R^{g+i}(\mathcal{L}_Q(M),C)=0\) for \(0<i<n-g.\)
\end{itemize}
\end{thm}

\pf:  Clear by Proposition \ref{GCSerreConditionAll} and Theorem \ref{horizontallylinkedtheorem}.\qed

To summarize some of these results
\begin{cor}
Let \((R,\mathfrak{m},k)\) be a local Cohen-Macaulay ring of dimension \(d\) with canonical module \(\omega_R\) and \(M\) an \(R\)-module of grade \(g.\)  Suppose \(M\) is directly linked to an \(\omega_R\)-quasi-Gorenstein module by the module $Q$.  Then the following hold
\begin{itemize}
\item[\((i)\))] \(M\in\mathcal{G}_{\omega_R}^g\)
\item[\((ii)\)] \(M\) is horizontally linked by \(Q\)
\item[\((iii)\)] \(\depth_{R_p}(M_p)\geq \depth(R_p)-g\) for all \(p\in \text{Spec}(R)\)
\item[\((iv)\)] \(\text{H}_{\mathfrak{m}}^i(M\otimes_R\omega_R)=0\) for \(g<i<d.\)
\end{itemize}
\end{cor}

As a closing remark, we will relate this summary to ideal linkage.  An important notion in ideal linkage is that of being in the linkage class of a complete intersection (licci).  This is analogous to being in the linkage class of an $R$-quasi-Gorenstein module in module linkage.  So the result above gives information  not only about modules which are linked through the canonical module $\omega_R$ but also about licci ideals. 

\vspace{0.5in}
\bibliographystyle{plain}

\begin{thebibliography}{10}

\bibitem{AuslanderBridgerSMT}
Maurice Auslander and Mark Bridger.
\newblock {\em Stable {M}odule {T}heory}.
\newblock Memoirs of the {A}merican {M}athematical {S}ociety, No. 94.
  {A}merican {M}athematical {S}ociety, {P}rovidence, {R}{I}, 1969.

\bibitem{RingHomoGDimAvramovFoxby}
Luchezar~L. Avramov and Hans-B{j\o\!\!}rn Foxby.
\newblock Ring homomorphisms and finite {G}orenstein dimension.
\newblock {\em Proc. London Math. Soc. (3)}, 75(2):241--270, 1997.

\bibitem{CIDIM}
Luchezar~L. Avramov, Vesselin~N. Gasharov, and Irena~V. Peeva.
\newblock Complete intersection dimension.
\newblock {\em Inst. Hautes \'Etudes Sci. Publ. Math.}, 86:67--114, 1997.

\bibitem{BrennanYorkQGorenstein}
Joseph Brennan and Alexander York.
\newblock quasi-{G}orensetin modules.
\newblock Thesis, 2017.

\bibitem{BrennanYorkETFS}
Joseph~P. Brennan and Alexander York.
\newblock An extension of a theorem of {F}robenius and {S}tickelberger to
  modules of projective dimension one over a factorial domain.
\newblock {\em J. Pure Appl. Algebra}, 223(2):626--633, 2019.

\bibitem{BH}
Winfried Bruns and J\"urgen Herzog.
\newblock {\em Cohen-{M}acaulay rings}, volume~39 of {\em Cambridge Studies in
  Advanced Mathematics}.
\newblock Cambridge University Press, Cambridge, 1993.

\bibitem{GDIMChristensen}
Lars~Winther Christensen.
\newblock {\em Gorenstein dimensions}, volume 1747 of {\em Lecture Notes in
  Mathematics}.
\newblock Springer-Verlag, Berlin, 2000.

\bibitem{CMDIMChristensenFoxbyFrankid}
Lars~Winther Christensen, Hans-Bj\o\!\!~rn Foxby, and Anders Frankild.
\newblock Restricted homological dimensions and {C}ohen-{M}acaulayness.
\newblock {\em J. Algebra}, 251(1):479--502, 2002.

\bibitem{GDimIyengar}
Lars~Winther Christensen and Srikanth Iyengar.
\newblock Gorenstein dimension of modules over homomorphisms.
\newblock {\em J. Pure Appl. Algebra}, 208(1):177--188, 2007.

\bibitem{LinkageModulesSerreSadDib}
Mohammad~T. Dibaei and Arash Sadeghi.
\newblock Linkage of modules and the {S}erre conditions.
\newblock {\em J. Pure Appl. Algebra}, 219(10):4458--4478, 2015.

\bibitem{EisenbudCA}
David Eisenbud.
\newblock {\em Commutative Algebra}.
\newblock Springer-Verlag, New York, 1995.

\bibitem{GIPFDIMEnochsJenda}
Edgar~E. Enochs and Overtoun M.~G. Jenda.
\newblock Gorenstein injective and flat dimensions.
\newblock {\em Math. Japon.}, 44(2):261--268, 1996.

\bibitem{AlgebraIIFaith}
Carl Faith.
\newblock {\em Algebra. {II}}.
\newblock Springer-Verlag, Berlin-New York, 1976.
\newblock Ring theory, Grundlehren der Mathematischen Wissenschaften, No. 191.

\bibitem{FossumDGR}
Robert Fossum.
\newblock Duality over {G}orenstein rings.
\newblock {\em Math. Scand.}, 26:165--176, 1970.

\bibitem{GorensteinModandRelatedFoxby}
Hans-Bj{\o\!\!}rn Foxby.
\newblock Gorenstein modules and related modules.
\newblock {\em Math. Scand.}, 31:267--284 (1973), 1972.

\bibitem{QuasiPerfectModFoxby}
Hans-Bj{\o\!\!}rn Foxby.
\newblock Quasi-perfect modules over {C}ohen-{M}acaulay rings.
\newblock {\em Math. Nachr.}, 66:103--110, 1975.

\bibitem{GDimandGenPerIGolod}
Evgenii~S. Golod.
\newblock {$G$}-dimension and generalized perfect ideals.
\newblock {\em Trudy Mat. Inst. Steklov.}, 165:62--66, 1984.
\newblock Algebraic geometry and its applications.

\bibitem{LiaisonHartshorne}
Robin Hartshorne.
\newblock Liaison with {C}ohen-{M}acaulay modules.
\newblock {\em Rend. Semin. Mat. Univ. Politec. Torino}, 64(4):419--432, 2006.

\bibitem{CMDIMHolmJorgensen}
Henrik Holm and Peter J\o\!\!~rgensen.
\newblock Cohen-{M}acaulay homological dimensions.
\newblock {\em Rend. Semin. Mat. Univ. Padova}, 117:87--112, 2007.

\bibitem{GenSerreHolmes}
Brent Holmes.
\newblock {A {G}eneralized {S}erre's {C}ondition}.
\newblock arXiv:1710.02631v1, October 2017.

\bibitem{StructureLinkageHuneke}
Craig Huneke and Bernd Ulrich.
\newblock The structure of linkage.
\newblock {\em Ann. of Math. (2)}, 126(2):277--334, 1987.

\bibitem{AlgebraicLinkageHuneke}
Craig Huneke and Bernd Ulrich.
\newblock Algebraic linkage.
\newblock {\em Duke Math. J.}, 56(3):415--429, 1988.

\bibitem{BasicAlgebra1}
Nathan Jacobson.
\newblock {\em Basic algebra. {I}}.
\newblock W. H. Freeman and Co., San Francisco, Calif., 1974.

\bibitem{BasicAlgebra2}
Nathan Jacobson.
\newblock {\em Basic algebra. {II}}.
\newblock W. H. Freeman and Co., San Francisco, Calif., 1980.

\bibitem{PureLinkageKustinMiller}
Andrew~R. Kustin, Matthew Miller, and Bernd Ulrich.
\newblock Linkage theory for algebras with pure resolutions.
\newblock {\em J. Algebra}, 102(1):199--228, 1986.

\bibitem{LinkageMartin}
Heath~M. Martin.
\newblock Linkage and the generic homology of modules.
\newblock {\em Comm. Algebra}, 28(9):4285--4301, 2000.

\bibitem{LinkageMartsinkovskyStrooker}
Alex Martsinkovsky and Jan~R. Strooker.
\newblock Linkage of modules.
\newblock {\em J. Algebra}, 271(2):587--626, 2004.

\bibitem{NagelLiaison}
Uwe Nagel.
\newblock Liaison classes of modules.
\newblock {\em J. Algebra}, 284(1):236--272, 2005.

\bibitem{SturgeonQGorenstein}
Uwe Nagel and Stephen Sturgeon.
\newblock Combinatorial interpretations of some {B}oij-{S}\"oderberg
  decompositions.
\newblock {\em J. Algebra}, 381:54--72, 2013.

\bibitem{LiaisonPeskineSzpiro}
Christian Peskine and Lucien Szpiro.
\newblock Liaison des vari\'et\'es alg\'ebriques. {I}.
\newblock {\em Invent. Math.}, 26:271--302, 1974.

\bibitem{SatherWagstaffNotes}
Sean Sather-Wagstaff.
\newblock Semidualizing modules.
\newblock unpublished NDSU course notes, 2009.

\bibitem{RelativeSemidualizingSather-Wagstaff}
Sean Sather-Wagstaff, Tirdad Sharif, and Diana White.
\newblock Comparison of relative cohomology theories with respect to
  semidualizing modules.
\newblock {\em Math. Z.}, 264(3):571--600, 2010.

\bibitem{SemidualizingSather-Wagstaff}
Sean Sather-Wagstaff and Siamak Yassemi.
\newblock Modules of finite homological dimension with respect to a
  semidualizing module.
\newblock {\em Arch. Math. (Basel)}, 93(2):111--121, 2009.

\bibitem{NotesLiaisonSchenzel}
Peter Schenzel.
\newblock Notes on liaison and duality.
\newblock {\em J. Math. Kyoto Univ.}, 22(3):485--498, 1982/83.

\end{thebibliography}

\end{document}